\documentclass[12pt, reqno]{amsart}
\usepackage{amscd,amsthm,amsfonts,amssymb,amsmath,euscript}
\usepackage[dvips]{graphicx}
\usepackage[dvips]{graphics}
\usepackage[matrix,arrow]{xy}

\newcounter{statements}

\theoremstyle{definition}

\theoremstyle{remark}
\newtheorem{remark}[statements]{Remark}

\newcommand{\QQ}{{\mathbb Q}}
\newcommand{\ZZ}{{\mathbb Z}}

\newcommand{\PP}{{\mathbb P}}
\newcommand{\CC}{{\mathbb C}}
\newcommand{\RR}{{\mathbb R}}

\newcommand{\Th}{{\bf Theorem }}
\newcommand{\Lem}{{\bf Lemma }}
\newcommand{\Def}{{\bf Definition }}
\newcommand{\Rem}{{\bf Remark}}

\newcommand{\Prop}{{\bf Proposition}}
\newcommand{\Proof}{{\bf Proof}}

\newfont{\smallskob}{cmbx7 scaled\magstep4}
\newfont{\bigskob}{cmbx12 scaled\magstep4}

\newcommand{\os}{\textup{\smallskob (}}
\newcommand{\cs}{\textup{\smallskob )}}

\newcommand{\qos}{\textup{\smallskob [}}
\newcommand{\qcs}{\textup{\smallskob ]}}

\newcommand{\itc}[1]{\textup{#1}}

\newcommand{\rk}{\mathrm{rk}\,}
\newcommand{\pic}{\mathrm{Pic}\,}

\newcommand{\virt}{\mathrm{virt}}
\newcommand{\tit}{Quantum cohomology of smooth complete intersections\\
in weighted projective spaces and singular toric varieties}

\begin{document}

\begin{title}
\tit
\end{title}

\begin{abstract}
We generalize Givental's Theorem for complete intersections in
smooth toric varieties in the Fano case. In particular, we find
Gromov--Witten invariants of Fano varieties of dimension $\geq 3$,
which are complete intersections in weighted projective spaces and
singular toric varieties. We generalize the Riemann--Roch
equations to
such varieties. We compute counting matrices of smooth Fano
threefolds with Picard group $\ZZ$ and anticanonical degrees 2, 8,
and 16.
\end{abstract}

\author{Victor Przyjalkowski}

\thanks{The work was partially supported by RFFI grants
$04-01-00613$ and $05-01-00353$, and grant NSh$-9969.2006.1$.}


\address{Steklov Institute of Mathematics, 8 Gubkin street, Moscow 119991, Russia} %




\email{victorprz@mi.ras.ru, przhijal@mccme.ru}

\maketitle

Let $X$ be a general complete intersection of hypersurfaces of
degrees $d_1,\ldots, d_l$ in $\PP^N$, $d_0=N+1-\sum d_i>0$,
$n=N-l\geqslant 3$. Let $H_H^*(X)\subset H^*(X,\ZZ)$ be the ring
generated by class $H$ that is dual to the hyperplane section and
$S=H_H^*(X)\otimes \CC[[q]]$. Denote $H\otimes 1\in S$, $H\in
H_H^*(X)$ by the same letter $H$. According to Givental's Theorem
for complete intersections in projective spaces
(see~\cite{Gi96}
), $I$-series for $X$ (i.e. the
generating series for one-pointed Gromov--Witten invariants with
descendants, see definition~\ref{definition:I-series}) is
$$
    {I^X}=e^{-\alpha_Xq}\sum_{d=0}^{\infty} \widehat I_d^X q^d=e^{-\alpha_Xq}\sum_{d=0}^{\infty}
    q^d \frac{\prod_{i=1}^l\prod_{j=1}^{dd_i} (d_iH+j)}{\prod_{i=1}^d
    (H+i)^{N+1}}\in S,
$$
where $\alpha_X=d_1!\ldots d_l!$ if $d_0=1$ and $\alpha_X=0$ if
$d_0>1$. Consider $D=q\cdot d/dq\in\mathcal D=\CC[[q]][q\cdot
d/dq]$. Then a solution of the Riemann--Roch equation (or equation
of type $Dn$), i. e. the differential equation
$$
(D^{N+1}-\prod_{i=0}^l d_i^{d_i}q \prod_{j=1}^{d_i}
(D+\frac{j}{d_i}))[\Phi(q)]=0
$$
is the series ${\mathcal I^X}=\sum q^d (d_0d)! \widehat I^X_{d,H^0}$. 

In this paper we generalize these results to complete
intersections in \emph{weighted projective spaces} and
\emph{singular toric varieties}.

In the first section we formulate the main results. In the second
section we give definitions and formulate known results. In the
third section we prove and discuss the main
theorems~\ref{theorem:main} and~\ref{theorem:toric}. In Appendix
we formulate Golyshev's conjecture. The corollary of
theorem~\ref{theorem:main} completes its proof.

\medskip

{\bf Notations.}  \emph{The Pochhammer symbol} $\os
X\cs_{\boldsymbol{n}}$ denotes the product
$X(X+1)\cdot \ldots \cdot (X+n-1)$, where $X$ belongs to some ring (see~\cite{AS72}, 6.1.22
). An infinite product of the form $\prod_{a=-\infty}^{n-1} (X+a)$
is denoted by $\qos X\qcs_{\boldsymbol{n}}$ (we define the
infinite product formally; we will only use quotients of two such
products 
in which all but for a finite number of factors coincide).

We denote the groups of homological and cohomological classes on
$X$ modulo torsion by $H_*(X)$ and $H^*(X)$ accordingly.

We use the same notation for a hypersurface and its dual
cohomological class.

Denote $B=\CC[[q]]$ and $D=q\cdot d/dq\in \mathcal D=B[q\cdot
d/dq]$.

\medskip

Everything is over $\CC$.

\section{Results}

See \ref{section:definition1}, \ref{section:definition2} and
\ref{section:definition3} for necessary definitions.

\subsection \Th
\label{theorem:main} {\it Let $\PP=\PP(w_1,\ldots,w_k)$ be a
weighted projective space. Let $X$ be a smooth complete
intersection of hypersurfaces $X_1,\ldots, X_l$ which do not
intersect the singular locus of $\PP$. Assume that $-K_X
> 0$ and $\pic X=\ZZ H$, where $H$ is the class dual to the hyperplane.
Let $i\colon X\rightarrow \PP$ be the natural embedding.

{\bf 1)} $I$-series
for $X$ is the
following.
$$
{I^X}=e^{-\alpha_X q}\sum_{d=0}^\infty q^{d}\cdot i^*\left(
\frac{\prod_{a=1}^l \os X_a+1 \cs_{\boldsymbol{d\cdot \deg X_a}}
}{\prod_{a=1}^k \os w_a H+1 \cs_{\boldsymbol{d\cdot w_a}}
}\right).
$$
Here $\alpha_X=0$ if the index of $X$ is $2$ or greater and
$\alpha_X=\prod_{a=1}^l (\deg X_a)!/\prod_{a=1}^k w_a!$ if the
index is $1$.

{\bf 2)} Let $d_i$ \itc{(}$1\leq i\leq l$\itc{)} be the degrees
\itc{(}with respect to $H$
\itc{)} of the hypersurfaces $X_i$ and $d_0=\sum w_i-\sum d_i$ be
the index of $X$.
Let ${\mathcal I^X}(q)=\sum q^d\cdot ((d_0d)!\cdot
(e^{\alpha_X}I^X)_{d,H^0})$.
Consider the operator \itc{(}which generalizes the Riemann--Roch
operator, see~\cite{Go01}\itc{)}
$$
L=\prod_{i=1}^{k} \os w_iD-(w_i-1)\cs_{\boldsymbol{w_i}} -q
\prod_{i=0}^{l} \os d_iD+1\cs_{\boldsymbol{d_i}}.
$$
Then $L[
\mathcal I^X(q)]=0$.}
\subsection \Th {\it

\label{theorem:toric} Let $Y$ be a $\QQ$-factorial toric variety
and $Y_1, \ldots, Y_k$ be the divisors that correspond to the
edges of the fan of $Y$. Consider a smooth complete intersection
$X$ of hypersurfaces $X_1,\ldots, X_l$ that does not intersect the
singular locus of $Y$. Assume that $-K_X
> 0$ and $\pic X=\ZZ$. Let $i\colon
X\rightarrow Y$ be the natural embedding. Let 
$\ell$ be a nef generator of $H_2(Y)$. For $\beta=d\ell$ put
$q^{\beta}=q^{d}$. Let $\Lambda\subset H_2(X)$ be the semigroup of
algebraic curves as cycles on $X$.

{\bf 1)} $I$-series
of $X$ is the
following.
$$
{I^X}=e^{-\alpha_X q}\sum_{\beta\in \Lambda} q^\beta\cdot
i^*\left( \frac{\prod_{a=1}^l \os X_a+1\cs_{\boldsymbol{\beta\cdot
X_a}} }{\prod_{a=1}^k \os Y_a+1 \cs_{\boldsymbol{\beta\cdot Y_a}}
}\right),
$$
where $\alpha_X=0$ if the index of $X$ is greater than $1$, and
$\alpha_X=\prod_{a=1}^l (\ell\cdot X_a)!/\prod_{a=1}^k (\ell\cdot
Y_k)!$ if the index is $1$.

{\bf 2)} Let $d_i$ \itc{(}$1\leq i\leq l$\itc{)} be the degrees of
the hypersurfaces $X_i$ \itc{(}with respect to $\ell$\itc{)}, $w_i$
\itc{(}$1\leq i\leq k$\itc{)} be the degrees of divisors that
correspond to the edges of the fan of $Y$, and $d_0$ be the index of
$X$.
Let ${\mathcal I^X}(q)=\sum q^d\cdot ((d_0d)!\cdot
(e^{\alpha_X}I^X)_{d,H^0})$.
Consider the operator
\itc{(}which generalizes the Riemann--Roch operator,
see~\cite{Go01}\itc{)}
$$
L=\prod_{i=1}^{k} \os w_iD-(w_i-1)\cs_{\boldsymbol{w_i}} -q
\prod_{i=0}^{l} \os d_iD+1\cs_{\boldsymbol{d_i}}.
$$
Then $L[
{\mathcal I^X}(q)]=0$.}

\begin{remark}
\label{remark:Lefschetz_condition} By Lefschetz Theorem
(see~\cite{Do82}, Theorem 3.2.4, (i) and Remark 3.2.6) hypotheses
of theorems~\ref{theorem:main} and~\ref{theorem:toric} hold for
complete intersections of dimension greater than $2$, which do not
intersect the singular locus of ambient variety, if its Picard
group is $\ZZ$ (this holds automatically by hypothesis of
theorem~\ref{theorem:main}).
\end{remark}

\section{Preliminaries}

Throughout the paper we consider the invariants of \emph{genus
zero} (those that correspond to the rational curves).

An axiomatic treatment of \emph{prime invariants} was given by
M.\,Kontsevich and Yu.\,Manin in~\cite{KM94}. \emph{Invariants with
descendants} were introduced and constructed in~\cite{BM96}
and~\cite{Beh97}.

\subsection{Moduli spaces of curves}
\label{section:definition1}
Consider smooth variety $X$ such that
$-K_X\geq 0$.

\subsubsection \Def {\it
\emph{The curve} is reduced scheme of pure dimension $1$. \emph{The
genus} of curve $C$ is the number $h^1(\mathcal O_C)$. }

It is easy to see that the curve is of genus $0$ if and only if it
is a tree of rational curves.

\subsubsection \Def {\it The connected curve $C$ with $n\geq 0$ marked
points $p_1,\ldots,p_n\in C$ is called \emph{prestable} if it has at
most ordinary double points as singularities and $p_1,\ldots, p_n$
are distinct smooth points \itc{(}see~\cite{Ma02}, III--2.1\itc{)}.
The map $f\colon C\rightarrow X$ of connected curve of genus $0$
with $n$ marked points are called \emph{stable} if $C$ is prestable
and there are at least three marked or singular points on every
contracted component of $C$ \itc{(}\cite{Ma02}, V--1.3.2\itc{)}. }

In the other words, a stable map of connected curve is the map that
has only finite number of infinitesimal automorphisms.

\subsubsection \Def {\it
\emph{The family of stable maps} \itc{(}over the scheme $S$\itc{)}
of curves of genus $0$ with $n$ marked points is the collection
$(\pi \colon \mathcal C\rightarrow S, p_1,\ldots,p_n, f\colon
\mathcal C\rightarrow X)$, where $\pi$ is the following map. It is a
smooth projective map with $n$ sections $p_1,\ldots,p_n$. Its
geometric fibers $(\mathcal C_s,p_1(s),\ldots,p_n(s))$ are prestable
curves of genus $0$ with $n$ marked points. Finally, the restriction
$f|_{\mathcal C_s}$ on each fiber is a stable map.

Two families over $S$
$$
(\pi\colon \mathcal C\rightarrow S, p_1,\ldots,p_n, f),\ \
(\pi^\prime\colon \mathcal C^\prime\rightarrow S,
p_1^\prime,\ldots,p_n^\prime, f^\prime)
$$
are called \emph{isomorphic} if there is an isomorphism $\tau\colon
\mathcal C\rightarrow \mathcal C^\prime$ such that
$\pi=\pi^\prime\circ \tau$, $p_i^\prime=\tau\circ p_i$,
$f=f^\prime\circ\tau$.
}

\subsubsection{} Let $\beta\in H_2^+(X)$. Consider the following
(contravariant) functor $ \overline{\mathcal M}_n(X, \beta)$ from
the category of (complex algebraic) schemes to the category of sets.
Let $\overline{\mathcal M}_n(X, \beta)(S)$ be the set of isomorphism
classes of families of stable maps of genus $0$ curves with $n$
marked points $(\pi\colon \mathcal C\rightarrow S, p_1,\ldots, p_n,
f)$ such that $f_*([\mathcal C_s])=\beta$, where $[C_s]$ is the
fundamental class of $C_s$.

\subsubsection \Def
{\it \emph{The moduli space of stable maps} of genus $0$ curves of
class $\beta\in H^+_2(X)$ with $n$ marked points to $X$ is the
Deligne--Mumford stack \itc{(}see~\cite{Ma02}, V--5.5\itc{)} which
is the coarse moduli space \itc{(}see~\cite{HM98}, Definition
$1.3$\itc{)} that represents $\overline{\mathcal M}_n(X, \beta)$.
This space is denoted by $\bar{M}_n(X,\beta)$. }


\subsection{Gromov--Witten invariants and $I$-series}
\label{section:definition2}

The definition of Gromov--Witten invariants is given in terms of
intersection theory on smooth stacks $\bar{M}_n(X,\beta)$. It exists
because locally such stack is a quotient of smooth variety by a
finite group (see~\cite{Vi89}). However such stacks may have
unexpected dimension. To use the products of cohomological cycles on
them one should introduce \emph{the virtual fundamental class}
$[\bar{M}_n(X,\beta)]^\virt$ of virtual dimension $\mathrm{vdim}\,
\bar{M}_n(X,\beta)=\dim X - \deg_{K_X} \beta+n-3$ (see its
construction in~\cite{Ma02}, VI--1.1). In many cases (for instance,
for homogenous spaces) virtual fundamental class coincides with the
usual one.

Consider the map $ev_i\colon \bar{M}_n(X,\beta)\rightarrow X$,
$ev_i(C;p_1,\ldots,p_n,f)=f(p_i)$. Let $\pi_{n+1}\colon
\bar{M}_{n+1}(X,\beta)\rightarrow \bar{M}_n(X,\beta)$ be the
forgetful map at the point $p_{n+1}$. It contracts unstable
components after forgetting $p_{n+1}$. Let $\sigma_i\colon
\bar{M}_n(X,\beta)\rightarrow \bar{M}_{n+1}(X,\beta)$ be the section
that are correspond to the marked point $p_i$ constructed as
follows. The image of the curve $(C;p_1,\ldots,p_n,f)$ under the map
$\sigma_i$ is the curve $(C';p_1,\ldots,p_{n+1},f')$. Here
$C'=C\bigcup C_0$, $C_0\backsimeq \PP^1$, and $C_0$ and $C$
intersect at the (non-marked on $C'$) point $p_i$. The points
$p_{n+1}$ and new $p_i$ lie on $C_0$. The map $f'$ contracts $C_0$
and $f'|_C=f$.

\begin{figure}[h]
  \includegraphics[scale=1]{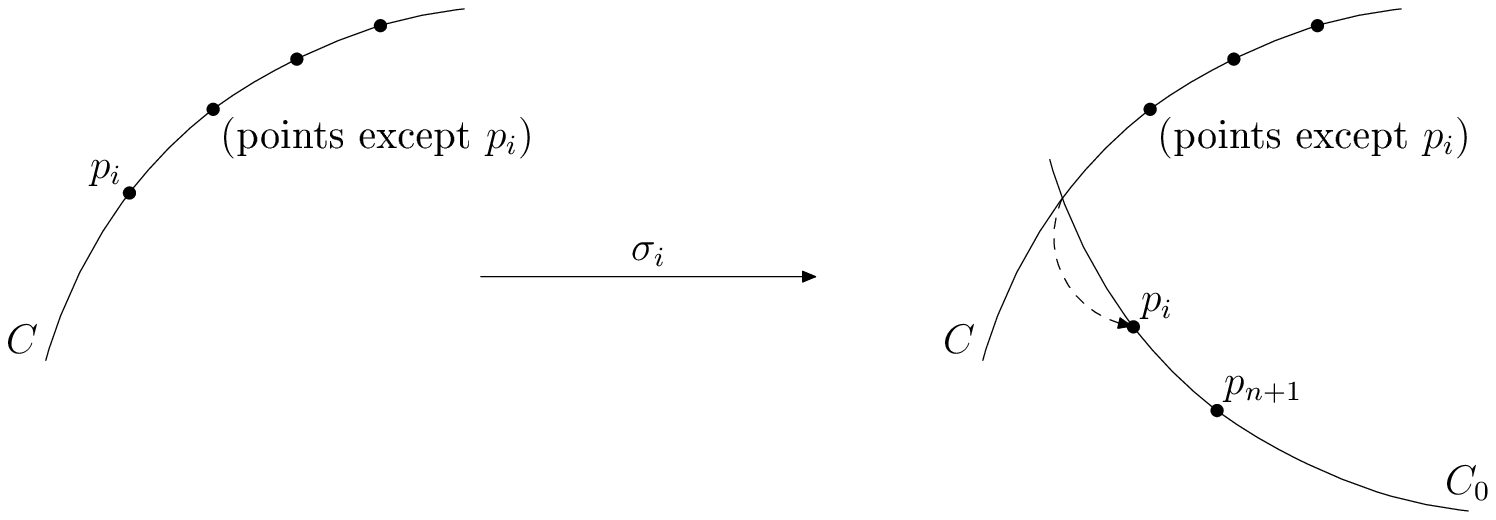}
  \label{fig}
\end{figure}

Let $L_i=\sigma_i^*\omega_{\pi_{n+1}}$, where $\omega_{\pi_{n+1}}$
is the relative dualizing sheaf $\pi_{n+1}$. The fiber $L_i$ over
the point $(C;p_1,\ldots,p_n,f)$ is $T^*_{p_i}C$.

\subsubsection{} \Def [see~\cite{Ma02}, VI--2.1]. {\it
\emph{A cotangent line class} is the class
$$
\psi_i=c_1(L_i)\in H^2(\bar{M}_n(X,\beta)).
$$ }

\subsubsection{} \Def [see \cite{Ma02}, VI--2.1]. {\it
\label{definition:GW-invariant} Consider the cohomological classes
$\gamma_1, \ldots, \gamma_n\in H^*(X)$. Let $a_1,\ldots, a_n$ be
the non-negative integers and $\beta \in H_2(X)$. \emph{The
Gromov--Witten invariant with descendants} that correspond to
these classes is
$$
\langle \tau_{a_1}\gamma_1, \ldots, \tau_{a_n}\gamma_n
\rangle_\beta= \psi_1^{a_1}ev_1^*(\gamma_1)\cdot\ldots\cdot
\psi_n^{a_n}ev_n^*(\gamma_n)\cdot [\bar{M}_n(X,\beta)]^\virt
$$
if $\sum \mathrm{codim}\, \gamma_i+\sum a_i=\mathrm{vdim}\,
\bar{M}_n(X,\beta)$ and $0$ otherwise. The invariants with $a_i=0$
\itc{(}for each $i$\itc{)} are called \emph{prime}. They are equal
to the expected numbers of rational curves of the class $\beta$ on
$X$, which intersect the dual cycles to $\gamma_1,\ldots,\gamma_n$
\itc{(}divided by the product of degrees of divisors from
$\gamma_1,\ldots,\gamma_n$\itc{)}. We omit the symbols $\tau_0$. }


\subsubsection{} \Def [see \cite{Ga00}
]. {\it \label{definition:I-series} Let $\mu_1,\ldots, \mu_N$ be the
basis of $H^*(X)$, $\check{\mu}_1,\ldots,\check{\mu}_N$ be the dual
basis, and $\gamma_1,\ldots, \gamma_k$ be the basis of $H_2(X)$.
Each curve $\beta$ is of the form $\beta=\sum \beta_i \gamma_i$.
Consider the ring $B=\CC[[q]]$, $q=(q_1,\ldots,q_k)$. Put
$q^\beta=\prod q_i^{\beta_i}$. Then \emph{the $I$-series} for $X$ is
the following.
\begin{gather*}
I^X=\sum_{\beta\geq 0} I_\beta^X\cdot q^\beta\in B\otimes H^*(X), \
\ \ \
I^X_\beta=
\sum_{i,j\geq 0} \langle\tau_i \mu_j\rangle_\beta\check{\mu}_j
\end{gather*}
\itc{(}we use the same symbol for $\gamma\in H^*(X)$ and
$1\otimes\gamma\in B\otimes H^*(X)$\itc{)}.

\emph{Fundamental term} of $I$-series
is the series
\begin{gather*}
I^X_{H^0}=\sum_{\beta\geq 0} \langle\tau_{(-K_X)\cdot\beta-2}
\check{\mathbf{1}}\rangle_\beta \cdot q^\beta,
\end{gather*}
where $\check{\mathbf{1}}$ is the class dual to the class of unity
in the cohomology.

}

\medskip

In the case $k=1$ and arbitrary series $I\in B\otimes H^*(X)$ put
$I=\sum_{q\geq 0} I_d \cdot q^d$.

\subsubsection{} \Lem [\cite{Ga99}, Lemma 5.5 or proof of Lemma 1
in~\cite{LP01}]. {\it \label{lemma:non-algebraic} Let $Y\subset X$
be a complete intersection and $\varphi \colon H^*(X) \rightarrow
H^*(Y)$ be the restriction homomorphism. Let
$\widetilde{\gamma}_1\in \varphi (H^*(X))^\bot$ and
$\gamma_2,\ldots, \gamma_l\in \varphi (H^*(X))$. Then for each
$\beta \in \varphi (H_2(X))\subset H_2(Y)$ the Gromov--Witten
invariant on $Y$ of the form
$$
\langle \tau_{d_1} \widetilde{\gamma}_1, \tau_{d_2} \gamma_2,
\ldots \tau_{d_l} \gamma_l \rangle_\beta
$$
vanishes. }

\subsubsection \Def {\it
Consider a complete intersection $Y \subset X$. Let $\varphi\colon
H^*(X)\rightarrow H^*(Y)$ be the restriction homomorphism and
$R=\varphi (H^*(X))\subset H^*(Y)$. Let $\mu_1,\ldots, \mu_N$ be the
basis of $R$, $\check{\mu}_1,\ldots,\check{\mu}_N\in R$ be the dual
basis, and $\gamma_1,\ldots, \gamma_k$ be the basis of $H_2(Y)$.
Each curve $\beta$ is of the form $\beta=\sum \beta_i \gamma_i$.
Consider the ring $B=\CC[[q]]$, $q=(q_1,\ldots,q_k)$. Put
$q^\beta=\prod q_i^{\beta_i}$.  Then \emph{restricted $I$-series} of
$Y$ is the following.
\begin{gather*}
\widetilde{I}^Y=\sum_{\beta\geq 0} \widetilde{I}_\beta^Y\cdot
q^\beta\in B\otimes R\subset B\otimes H^*(Y),\ \ \ \
\widetilde{I}^Y_\beta=
\sum_{i,j\geq 0} \langle\tau_i \mu_j\rangle_\beta\check{\mu}_j.
\end{gather*}

The invariants of complete intersection that correspond to
cohomological classes restricted from the ambient variety are
called restricted.
}

\subsubsection{} \Rem.
\label{remark:restricted_I-series} For a complete intersection
$Y\subset X$ of dimension at least $3$ these series coincide, i.~e.
$I^Y=\widetilde{I}^Y$. Indeed, $H_2(Y)\simeq H_2(X)$ and by
lemma~\ref{lemma:non-algebraic} and the divisor
axiom~\ref{axiom:divisor} one-pointed invariants for primitive
classes vanish.

\subsubsection \Rem
Two- and more- pointed invariants for primitive classes of middle
dimension may be non-zero. For instance, for two primitive classes
$\alpha$ and $\beta$, a hypersurface section $H$ and a line $\ell$
on cubic threefold
$$
\langle \alpha, \beta, H^2 \rangle_\ell=-\alpha\cdot\beta
$$
(see~\cite{Bea95}, Proposition 1).


\subsection{Relations between multipointed invariants}

\subsubsection{}
In this paragraph we discuss relations between one-pointed
Gromov--Witten invariants with descendants (i.~e. coefficients of
$I$-series) and multipointed prime ones (in particular, two-pointed
ones, which are important for the following).

One-pointed invariants with descendants may be expressed in terms of
multipointed ones by the following.

\subsubsection{} {\bf Divisor axiom for invariants with descendants} [\cite{KM98}, 1.5.2].
\label{axiom:divisor} {\it Let $Y$ be smooth projective
variety, $\gamma_1,\ldots, \gamma_n\in H^*(Y)$ and $\gamma_0\in
H^2(Y)$ be a divisor. Then
\begin{multline*}
\langle\gamma_0, \tau_{d_1} \gamma_1,\ldots,
\tau_{d_n}\gamma_n\rangle_\beta=
(\gamma_0\cdot\beta)\langle\tau_{d_1} \gamma_1,\ldots, \tau_{d_n}\gamma_n\rangle_\beta+ \\
\sum_{k,d_k\geq 1} \langle\tau_{d_1} \gamma_1,\ldots,
\tau_{d_k-1}(\gamma_0\cdot
\gamma_k),\ldots,\tau_{d_n}\gamma_n\rangle_\beta.
\end{multline*}
}

\medskip

The following formula (with divisor axiom) enables one to express
recursively three-pointed invariants with descendants in terms of
prime ones.

\subsubsection{} \Th [topological recursion relations,
\cite{Ma02}, VI--6.2.1]. \label{theorem: three-point-correlation}
{\it Let $Y$ be a smooth projective variety, $\{\Delta^i\}$ and
$\{\Delta_i\}$, $i=1,\ldots, N$, be dual bases of $H^*(Y)$. Then
$$
\langle \tau_{d_1} \gamma_1, \tau_{d_2} \gamma_2, \tau_{d_3}
\gamma_3 \rangle_\beta = \sum_{a, \beta_1+\beta_2=\beta} \langle
\tau_{d_1-1} \gamma_1, \Delta^a \rangle_{\beta_1} \langle \Delta_a,
\tau_{d_2} \gamma_2, \tau_{d_3} \gamma_3 \rangle_{\beta_2}
$$
\itc{(}the sum is taken over all $a,\beta_1,\beta_2$ such that the
expression in it is well-defined, i.~e. $a=1,\ldots,N$,
$\beta_1,\beta_2\geq 0$\itc{)}.

}

\subsubsection{} Using this formulas, we can express one-pointed
invariants in terms of prime two-pointed ones (see~\cite{Pr04},
Proposition $5.2$). In fact, in many cases these expressions are
invertible. More particular, the inverse expressions can be found
recursively by the divisor axiom and the following theorem (see
also~\cite{BK00}, теорема $5.2$ and~\cite{Pr04}).



\subsubsection{} \Th [Lee, Pandharipande,~\cite{LP01}, Theorem 
2,\,i]. \label{theorem:multipoint expressions} {\it Let $Y$ be a
smooth projective variety. Consider the self-dual ring $R\subset
H^*(Y)$ generated by Picard group $\pic Y$ such that for any
$\mu_1,\ldots,\mu_n\in R$ and $\nu\in R^\bot$
$$
\langle \tau^{a_1} \mu_1,\ldots, \tau^{a_n} \mu_n, \tau^a \nu
\rangle_\beta=0.
$$
Let $\gamma_1,\ldots,\gamma_n\in R$, $H\in \pic Y$, $\Delta^i$ and
$\Delta_i$ be the dual bases of $H^*(Y)$. Then one can algebraically
express Gromov--Witten invariant
$$
\langle \tau^{k_1} \gamma_1,\ldots, \tau^{k_n} \gamma_n
\rangle_\beta
$$
in terms of one-pointed Gromov--Witten invariants with descendants
using the following expressions
\begin{multline*}
\langle \tau_{k_1} \gamma_1,\ldots, \tau_{k_n} H\cdot\gamma_n
\rangle_\beta=\langle \tau_{k_1} H\cdot\gamma_1,\ldots, \tau_{k_n}
\gamma_n \rangle_\beta+\beta H\cdot \langle \tau_{k_1+1}
\gamma_1,\ldots, \tau_{k_n} \gamma_n \rangle_\beta-\\
\sum_{\beta_1+\beta_2=\beta} \beta_1 H\cdot \sum_{S^1\cup S^n=S,\
a } \langle \tau_{k_{s^1_1}} {\gamma_{s^1_1}},\ldots,
\tau_{k_{s^1_a}} \gamma_{s^1_a}, \Delta^a
\rangle_{\beta_1}\cdot\langle \Delta_a, \tau_{k_{s^n_1}}
{\gamma_{s^n_1}},\ldots, \tau_{k_{s^n_b}} \gamma_{s^n_b}
\rangle_{\beta_2}
\end{multline*}
and
\begin{multline*}
\langle \tau_{k_1} \gamma_1,\ldots, \tau_{k_n+1} \gamma_n
\rangle_\beta=-\langle \tau_{k_1+1} \gamma_1,\ldots,
\tau_{k_n} \gamma_n \rangle_\beta+\\
\sum_{
\parbox[c]{1,2cm}{\scriptsize
$\beta_1+\beta_2=\beta,$ $S^1\cup S^n=S,\ a$}} \langle
\tau_{k_{s^1_1}} {\gamma_{s^1_1}},\ldots, \tau_{k_{s^1_a}}
\gamma_{s^1_a}, \Delta^a \rangle_{\beta_1}\cdot\langle \Delta_a,
\tau_{k_{s^n_1}} {\gamma_{s^n_1}},\ldots, \tau_{k_{s^n_b}}
\gamma_{s^n_b} \rangle_{\beta_2},
\end{multline*}
where $H\in \pic X\cap R$ and the latter summations are taken over
partitions $S^1=\{s^1_1,\ldots, s^1_a\}$ and $S^n=\{s^n_1,\ldots,
s^n_b\}$ of the set $S=\{1,\ldots,n\}$ such that $1\in S^1$ and
$n\in S^n$. }

\subsubsection{}


The hypothesis of this theorem holds for subring of restricted
cohomologies of complete intersection in toric variety
(see~\cite{Fu93}, 5.2) or for the subring of the cohomology ring
of Fano threefold generated by Picard group. See expressions for
Fano threefolds with Picard group $\ZZ$ in~\cite{Pr04}.

\emph{ Consequently, theorems~\ref{theorem:main}
and~\ref{theorem:toric} enable us to find the restricted quantum
cohomology ring of smooth complete intersections in weighted
projective spaces and singular toric varieties.}

\subsection{Toric varieties}
\label{section:definition3}

The definition and the main properties of toric varieties see
in~\cite{Da78} or in~\cite{Fu93}. Just remind that toric variety is
a variety with action of torus $T\simeq (\CC^*)^n$ such that one of
its orbits is a Zariski open set. Toric variety is determined by its
\emph{fan}, i. e. some collection of cones with vertices in the
points of lattice that is dual to the lattice of torus characters.
Moreover, algebraic-geometric properties of toric variety can be
formulated in terms of properties of this fan.

Remind some of them.

\subsubsection{}
Every cone of the fan $\sigma\subset N=\RR^n$ of dimension $r$
corresponds to the orbit of the torus of dimension $n-r$ (here $n$
is the dimension of toric variety). Thus, each edge (one-dimensional
cone) correspond to the (equivariant) Weil
divisor\footnote{\label{toric} That is, let $\Sigma\in N=\ZZ^n$ be a
fan of the toric variety $X_\Sigma$ and let $\sigma\in \Sigma$ be
any cone. Let $M$ be a lattice dual to $N$ with respect to some
non-degenerate pairing $\langle \cdot,\cdot\rangle$ and
$\sigma^\vee$ be a dual cone for $\sigma$ (i.~e. $\sigma^\vee=
\{l\in M| \forall k\in \sigma \ \ \langle l,k\rangle\geq 0\}$). Let
$U_\sigma=\mathrm{Spec}\, \CC [\sigma^\vee]$ corresponds to
$\sigma$. The variety $X_\sigma$ is obtained from the affine
varieties $U_\sigma$, $\sigma\in \Sigma$, by gluing together
$U_\sigma$ and $U_\tau$ along $U_{\sigma\cap\tau}$, $\sigma,\tau\in
\Sigma$. Thus, if $l\subset\sigma\in \Sigma$ is an edge of the fan,
then the divisor that is correspond to $l$ restricted on $U_\sigma$
as $U_l\subset U_\sigma$. }.
The divisors which correspond to the edges of the fan generate
divisor class group. A Weil divisor $D=\sum d_i M_i$, where $M_i$
corresponds to edges, is Cartier if for each cone of the fan
$\sigma$ there exist a vector $n_\sigma$ such that $\langle
n_\sigma, m_i\rangle=d_i$ where $m_i$ is the primitive elements of
the edges of this cone. If such vector is the same for all cones,
then the divisor is principal. Hence if the toric variety is
$n$-dimensional and the number of the edges is $k$, then the rank of
the divisor class group is $n-k$.

\subsubsection \Def {\it
The variety is called \emph{$\QQ$-factorial} if for each Weil
divisor $D$ there exist some integer $k$ such that $kD$ is a Cartier
divisor. }

\medskip

In particular, there exist an intersection theory for Weil divisors
on the $\mbox{$\QQ$-factorial}$ variety.

\subsubsection{} Toric variety is $\QQ$-factorial if and only if any cone of the
fan, which corresponds to this variety, is simplicial. In this
case Picard group is generated (over $\QQ$) by divisors, which
correspond to the edges of the fan.

\subsubsection{} Consider a weighted projective space
$\PP=\PP(w_0,\ldots,w_l)$. The fan which correspond to it is
generated by integer vectors $m_0,\ldots,m_l\in \RR^l$ such that
$\sum w_i m_i=0$. If $w_0=1$, then one can put $m_0=(-w_1,\ldots,
-w_l)$, $m_i=e_i$, where $e_i$ is a basis of $\RR^l$. The collection
$\{m_i\}$ corresponds to the collection of standard divisors--strata
$\{ D_i\in
|w_iH|
\}$.

\subsubsection{} A toric variety is \emph{smooth}
if for any cone $\sigma$ in the fan that correspond to this
variety the subgroup $\sigma\cap \ZZ^n$ is generated by the subset
of the basis of the lattice $m_1^\sigma, \ldots, m_k^\sigma$.
Adding the edge $a=a_1 m_1^\sigma+\ldots+a_km_k^\sigma$, $a_i\in
\QQ$ to the cone (and connection it with ``neighboring'' faces)
corresponds to weighted blow-up (along subvariety which correspond
to $\sigma$) with weights $1/r\cdot(\alpha_1,\ldots,\alpha_k)$,
where $\alpha_i\in \ZZ$ and $a_i=\alpha_i/r$. Consecutively adding
edges to the fan in this way we can get toric resolution of a
toric variety.

\subsubsection{}
\label{remark:singular-weighted-projective} So, singular locus of
weighted projective space $\PP=\PP(w_0,\ldots,w_l)$ is the union of
strata given by $x_{i_1}=\ldots=x_{i_j}=0$, where $x_{i_j}$ is the
coordinate of weight $w_{i_j}$ and $\{ i_1,\ldots,i_j \}$ is the
maximal set of indices such that greatest common factor of the
others is greater than $1$.

%
%
%
%
%
%
%

\subsection{} {\bf Givental's Theorem} [\cite{Gi97}, Theorem $0.1$].
{\it \label{theorem:givental} Let $X$ be a smooth toric variety
and $Y$ be a smooth complete intersection in it with positive
anticanonical class. Let $X_1,\ldots,X_k$ be the divisors which
correspond to the edges of a fan of $X$ and let $Y$ be given by
divisors $Y_1,\ldots, Y_r$. Let $\Lambda\subset H_2(X)$ be the
semigroup of algebraic curves as cycles on $Y$, and $i\colon
Y\rightarrow X$ be a natural embedding. Then
$$
I^Y=e^{h(q)}\sum_{\beta\in \Lambda}q^\beta \cdot i^* \left(
\frac{\prod_{a=1}^r \qos Y_a \qcs_{\boldsymbol{\beta\cdot Y_a+1}}
}{\prod_{a=1}^r \qos Y_a \qcs_{\boldsymbol{1}} }
\frac{\prod_{a=1}^k \qos X_a \qcs_{\boldsymbol{1}} }{\prod_{a=1}^k
\qos X_a \qcs_{\boldsymbol{\beta\cdot X_a+1}} }\right),
$$
where $h(q)$ is a polynomial supported by curves, whose
intersection with anticanonical class is $1$ \itc{(}i. e.
$h(q)=\sum_\beta h_\beta q^\beta$, where $\beta \cdot
(-K_Y)=1$\itc{)}. }

\subsection{Motivation}
Consider a smooth Fano threefold $X$ with Picard group $\ZZ$. Put
$K=-K_X$.

\subsubsection{} \Def [\cite{Go02}, 1.7, 1.10].
\label{definition:counting-matrix} {\it \emph{A counting matrix} is
the matrix of Gromov-Witten invariants of $X$, namely
the following matrix $A\in \mathrm{Mat}(4\times 4)$.
$$
    A=\left[\begin{array}{cccc}
    a_{00} & a_{01} & a_{02} & a_{03} \\
    1      & a_{11} & a_{12} & a_{13} \\
    0      & 1      & a_{22} & a_{23} \\
    0      & 0      & 1      & a_{33} \\
    \end{array}\right].
$$
Numeration of rows and columns starts from $0$ and the elements are
given by
$$
a_{ij}=\frac{\langle K^{3-i}, K^j, K\rangle_{j-i+1}}{\deg
X}=\frac{j-i+1}{\deg X}\cdot\langle K^{3-i}, K^j \rangle_{j-i+1}$$
\itc{(}the degree is taken with respect to the anticanonical
class\itc{)}. }

\subsubsection{} It is easy to see that the matrix $A$ is symmetric with
respect to the secondary diagonal: $a_{ij}=a_{3-j,3-i}$. By
definition, $a_{ij}=0$ if $j-i+1<0$. If $j-i+1=0$, then $a_{ij}=1$,
because it is just a number of intersection points of $K^{3-i}$,
$K^j$, and $K$, which is $\deg X$; $a_{00}=a_{33}=0$. For the other
coefficients $a_{ij}$'s are ``expected'' numbers of rational curves
of degree $j-i+1$ passing through $K^{3-i}$ and $K^j$, multiplied by
$\frac{j-i+1}{\deg X}$. The only exception is the following: by
divisor axiom
$$
a_{01}=2\cdot (2\cdot ind\,(X)\cdot\mathrm{[\emph{the number of
conics passing through the general point}])}.
$$

\subsubsection{} Consider the following Fano threefolds \itc{(}see~\cite{Is77}, \cite{Is78}, \cite{Is79},
\cite{Is88}, \cite{IP99}, \cite{Mu92}\itc{)}.

\begin{description}

\item[1]
The variety $V_1$ of anticanonical degree 8 \itc{(}a double
covering of the cone over Veronese surface branched over a smooth
cubic\itc{)}.

\item[2]
The variety $V_2$ of anticanonical degree 16 \itc{(}a double
covering of $\PP^3$ branched over a smooth quartic\itc{)}.

\item[3]
The variety $V_2'$ of anticanonical degree 2 \itc{(}a double
covering of $\PP^3$ branched over a smooth sextic\itc{)}.

\end{description}

\subsubsection{} \Prop. {\it
\label{proposition:hypersurfaces} The varieties $V_1$, $V_2$, and
$V_2'$ are of the following form.

\begin{description}

\item[1]
Any variety of type $V_1$ is isomorphic to a smooth hypersurface of
degree $6$ in $\PP(1,1,1,2,3)$.

\item[2]
Any variety of type $V_2$ is isomorphic to a smooth hypersurface of
degree $4$ in $\PP(1,1,1,1,2)$.

\item[3]
Any variety of type $V_2'$ is isomorphic to a smooth hypersurface of
degree $6$ in $\PP(1,1,1,1,3)$.

\end{description}
}

\subsubsection \Proof
Double covering $\PP(w_0,\ldots,w_n)$ branched over a divisor
given by function $f_k(x_0,\ldots,x_n)$ of degree $k$ may be given
by function $x_{n+1}^2=f_k$ in $\PP(w_0,\ldots,w_n,k/2)$. The
variables $x_i$ here have the weights $w_i$ and the variable
$x_{n+1}$ has the weight $k/2$. \qed

\subsubsection \Th {\it
\label{theorem:matrices} The counting matrices for $V_1$, $V_2$,
and $V_2'$ are as follows.

\begin{description}

\item[1]
  For $V_{1}$
$$
    M(V_{1})=
    \left[\begin{array}{cccc}
         0 &     240&     0   &   57600  \\
    1      &    0   &    1248 &     0    \\
    0      & 1      &     0   &     240  \\
    0      & 0      & 1       &     0    \\
    \end{array}\right].
$$

\item[2]
  For $V_{2}$
$$
    M(V_{2})=
    \left[\begin{array}{cccc}
         0 & 48 &  0 & 2304 \\
    1      &   0  &  160    & 0   \\
    0      & 1      &    0     & 48      \\
    0      & 0      & 1          &     0       \\
    \end{array}\right].
$$

\item[3]
  For $V_{2}'$
$$
    M(V_{2}')=
    \left[\begin{array}{cccc}
         0 &     137520 &     119681240 &    21690374400   \\
    1      &    624   &     650016  &     119681240  \\
    0      & 1      &     624   &     137520   \\
    0      & 0      & 1       &     0    \\
    \end{array}\right].
$$

\end{description}  }

\subsubsection{} {\bf Proof.}
By proposition \ref{proposition:hypersurfaces} these varieties are
hypersurfaces in weighted projective spaces. Find their one-pointed
invariants using theorem~\ref{theorem:main}. Apply
theorem~\ref{theorem:multipoint expressions} for the subring of
cohomology ring generated by the class dual to the class of
hyperplane section and find prime two-pointed invariants. They are
the coefficients of counting matrices. \qed


\section{Proofs of the main theorems}

\subsection{} {\bf Proof of theorems~\ref{theorem:main} and~\ref{theorem:toric}}.

The divisors in $\PP$ that correspond to the edges of its fan are
$w_1H,\ldots,w_kH$. Thus, theorem~\ref{theorem:main} follows from
theorem~\ref{theorem:toric}.

Prove theorem~\ref{theorem:toric}.

{\bf 1)} Let $\Sigma$ be a fan of $Y$. Consider the sequence of
simplicial fans $\Sigma=\Sigma_0, \Sigma_1,\ldots, \Sigma_r$ such
that
\begin{itemize}
  \item[i)] The fan $\Sigma_{i+1}$ is a result of the following procedure
applied to $\Sigma_i$. Pick a cone whose edges are not part of a
basis of the integer lattice containing this fan. Add an edge,
which lies inside the cone. Add the faces which contain this edge
and ``neighboring'' edges (i. e. replace all cones which contain
the edge we add by all linear spans of this edge and faces of this
cone that does not contain the edge);
 \item[ii)]
The fan $\Sigma_r$ corresponds to a smooth variety.
\end{itemize}
Such sequence exists by~\cite{Da78}, $8.1$--$8.3$. 

This way we get a toric resolution of singularities $f\colon
\widetilde{Y}=X_{\Sigma_r}\rightarrow Y$ (where $X_{\Sigma_r}$ is
the toric variety that correspond to the fan $\Sigma_r$, see
remark~\ref{toric}). The exceptional set of this resolution is the
union of divisors $E_1,\ldots,E_r$ which correspond to the edges we
add. Let $W_1,\ldots,W_k$ be divisors that correspond to the edges
of a fan of $Y$ and $\widetilde{W}_1,\ldots,\widetilde{W}_k$ be
their strict transforms. Let
$\widetilde{X}_1,\ldots,\widetilde{X}_l$ be strict transforms of
$X_1,\ldots,X_l$,
$\widetilde{X}=\widetilde{X}_1\cap\ldots\cap\widetilde{X}_l$, and
$j\colon \widetilde{X}\rightarrow \widetilde{Y}$ be the natural
embedding. Obviously, there exists an isomorphism $g$ such that the
diagram
\[ \xymatrix{
\widetilde{X}\ar@{->}[d]_{g}\ar@{->}[rr]^{j}&& \widetilde{Y}\ar[d]^{f}\\
X\ar@{->}[rr]^{i}&&Y}\] is commutative, because $f$ is an
isomorphism between $\widetilde{Y}\setminus (\cup_{i=1}^r E_i)$
and $Y\setminus \mathrm{Sing}\, Y$.

Consider any toric variety $V$ associated to a complete simplicial
fan. Then the canonical map of group of algebraic cycles modulo
rational equivalence with rational coefficients to homology group
with rational coefficients
$$
A_*(V)_{\QQ}\rightarrow H_*(V,\QQ)
$$
is an isomorphism (see~\cite{Da78}, 10.9). So, we can extend an
intersection theory to the homology group with rational
coefficients.

The map $f_*\colon H_2(\widetilde{Y},\QQ)\rightarrow H_2(Y,\QQ)$
is surjective. Let $\lambda\subset X\subset Y$ be an effective
curve such that $H_2(X)=\ZZ\lambda$ and
$\widetilde{\lambda}=f^{-1}(\lambda)$. Put $K=\ker f_*$. Then
$H_2(\widetilde{Y},\QQ)=\QQ\widetilde{\lambda}+K$.

The cycle of any curve $\boldsymbol{\beta}$ which lies on
$\widetilde{X}$ equals $\beta_0 \cdot \widetilde{\lambda}$, where
$\beta_0\in \ZZ$ and $\beta_0\geq 0$, because representatives of
elements of $K$ lie on exceptional divisors, which do not
intersect $\widetilde{X}$. So, if $\widetilde{\Lambda}$ is a
semigroup of algebraic curves on $\widetilde{X}$, then
$\widetilde{\Lambda}=\ZZ_{\geq 0}\widetilde{\lambda}$.

The multiplicity of intersection of divisor and a curve depends
only on the neighborhood of this curve and $f$ is an isomorphism
in a neighborhood of $\widetilde{X}$, so $\widetilde{\lambda}\cdot
\widetilde{X}_i=\lambda\cdot X_i$ and $\widetilde{\lambda}\cdot
\widetilde{W}_i=\lambda\cdot W_i$. Thus, for any curve
$\boldsymbol{\beta}=\beta_0\cdot \widetilde{\lambda}$ on
$\widetilde{X}$ we have $\boldsymbol{\beta}\cdot
\widetilde{W}_i=\beta_0\cdot (\lambda\cdot W_i)$,
$\boldsymbol{\beta}\cdot \widetilde{X}_i=\beta_0\cdot
(\lambda\cdot X_i)$, and $\boldsymbol{\beta}\cdot E_i=0$.

Divisors that correspond to the edges of a fan of $\widetilde{Y}$
are $\widetilde{W}_1,\ldots, \widetilde{W}_k, E_1,\ldots, E_r$. But
the factors in the expression for $I$-series for $\widetilde{X}$
from Theorem~\ref{theorem:givental} which include the divisors $E_i$
may be cancelled, so we may omit them. Consider rings $H^*(X)[[q]]$
and $H^*(\widetilde{X})[[\widetilde{q}]]$ and an isomorphism between
them which sends $q$ to $\widetilde{q}$ and acts on the coefficients
as $f^*$. Let $j\colon \widetilde{X}\rightarrow \widetilde{Y}$ be a
natural embedding. Then

\begin{multline*}
I^{\widetilde{X}}=
e^{\widetilde{h}(\widetilde{q})}\sum_{\beta=\beta_0\widetilde{\lambda}\in
\widetilde{\Lambda} } \widetilde{q}^{\beta_0} \cdot
j^*\left(\frac{\prod_{a=1}^l \qos \widetilde{X}_a
\qcs_{\boldsymbol{\beta\cdot \widetilde{X}_a+1}} }{\prod_{a=1}^l
\qos \widetilde{X}_a\qcs_{\boldsymbol1} } \frac{\prod_{a=1}^k \qos
\widetilde{W}_a\qcs_{\boldsymbol1} }{\prod_{a=1}^k \qos
\widetilde{W}_a\qcs_{\boldsymbol{\beta\cdot \widetilde{W}_a+1}
}}\right)
=\\
j^*f^{*}\left(e^{h(q)}\sum_{\beta_0\geq 0} q^{\beta_0} \cdot
\prod_{a=1}^l \os X_a+1\cs_{\boldsymbol{\beta_0\cdot \deg X_a}}
\cdot \frac{1}
{\prod_{a=1}^k \os \widetilde{W}_a+1\cs_{\boldsymbol{\beta_0\cdot w_a}} }\right)=\\
g^*\left(e^{h(q)}\sum_{d=0}^\infty q^{d} \cdot
i^*\left(\frac{\prod_{a=1}^l \os X_a+1\cs_{\boldsymbol{d\cdot \deg
X_a}} }{\prod_{a=1}^k \os \widetilde{W}_a+1\cs_{\boldsymbol{d\cdot
w_a}} }\right)\right),
\end{multline*}
where $\widetilde{h}(\widetilde{q})$ is a polynomial from
theorem~\ref{theorem:givental} and $h(q)$ is its pre-image. The
second equality holds because we may choose divisors that
intersect $X$ properly, and their strict transforms as
representatives (over $\QQ$) of divisor classes (because $\pic
Y=\ZZ$ and every effective divisor on $Y$ is ample, and,
therefore, it is equivalent to a rational multiple to a very ample
one). The map $f$ after restriction on $\widetilde{X}$ is
isomorphism, which induces the isomorphism $f^*$ of cohomology
groups such that the image of any divisor under $f^*$ is its
strict transform. Finally, $g^*$ is an isomorphism of cohomology
rings $X$ and $\widetilde{X}$, so $I^{\widetilde{X}}=g^*(I^X)$.

The only thing that has to be found is $h(q)$. If the index of $X$
is greater than $1$, then there are no curves which intersect the
anticanonical class by $1$, so $h(q)=0$. If the index is 1, then any
such curve is a line $\ell$ (with respect to the anticanonical class
of $X$) and $h(q)=-\alpha_X q$. Besides, there are no lines on $X$
passing through the general point, i. e. $\langle H^{\dim X}
\rangle_\ell=0$. This means that the coefficient $e^{-\alpha_X
q}\cdot I^X_\ell$ at $H^0$ vanishes, so
$$
\alpha_X=\frac{\prod_{a=1}^l (\ell\cdot X_a)!}{\prod_{a=1}^k
(\ell\cdot Y_k)!}=\frac{\prod_{a=1}^l (\deg X_a)!}{\prod_{a=1}^k
w_a!}.
$$

{\bf 2)} It is easy to see that
$$
\mathcal I^X(q)_{H^0}=\sum_{n=0}^\infty \frac{\prod_{i=0}^l
(d_i\cdot n)!}{\prod_{j=1}^k (w_j\cdot n)!} q^n.
$$
Consider the solution of the equation $L[\mathcal J^X(q)]=0$ as a
series $\mathcal J^X=\sum_{n\geq 0} a_nq^n$. Evidently,
$D[q^n]=nq^n$. Hence
$$
(\prod_{i=1}^{k} 
\os w_iD-(w_i-1)\cs_{\boldsymbol{w_i}} )[q^n]= (\prod_{i=1}^{k}
\os w_in-(w_i-1)\cs_{\boldsymbol{w_i}})\cdot q^n.
$$
Analogously,
$$
(q \prod_{i=0}^{l} 
\os d_iD+1\cs_{\boldsymbol{d_i}} )[q^{n-1}]=(\prod_{i=0}^{l} \os
d_i(n-1)+1\cs_{\boldsymbol{d_i}} )\cdot q^n.
$$
The differential equation on $\mathcal J^X$ translates into
recursive equations on the coefficients $a_n$ using these
equalities. Namely,
$$
a_n=\frac{\prod_{i=0}^{l} \os d_i(n-1)+1\cs_{\boldsymbol{d_i}}
}{\prod_{i=1}^{k} \os w_in-(w_i-1)\cs_{\boldsymbol{w_i}}} \cdot
a_{n-1}.
$$
Put $a_0=1$. By induction we get
$$
a_n=\sum \frac{\prod_{i=0}^l (d_i\cdot n)!}{\prod_{j=1}^k
(w_j\cdot n)!} q^n.
$$
\qed


\subsubsection \Rem

It is easy to see that the Riemann--Roch operator is divisible by
$D^{\dim Y-\dim X}$ from the left. Let $L=D^{\dim Y-\dim
X}\widetilde{L}$. Then $\widetilde{L}\mathcal I^X=0$.

\subsection{} {\bf Comments and examples}.

\subsubsection{} {\bf Remarks.}
\begin{description}
    \item[i]
In order to check that a complete intersection $X=X_1\cup \ldots
\cup X_l$ of general hypersurfaces of degrees $d_1,\ldots, d_l$ in a
weighted projective space $\PP=\PP (w_1, \ldots w_k)$ does not
intersect the singular locus of $\PP$ one can use the following
necessary (but not sufficient) condition.

\begin{multline*}
\mbox{\it The number of Cartier divisors among $X_i$ is greater
than} \\
\mbox{\it the dimension of the singular locus of $\PP$.}
\end{multline*}
This means that

\begin{multline*}
    (\mbox{\it a number of $d_i$'s that is divisible on each of $a_i$'s})> (\mbox{\it the maximal number of weights
    $w_i$}\\
    \mbox{\it such that greatest common factor of the others is greater than 1})
\end{multline*}
(see~\ref{remark:singular-weighted-projective}).

For an arbitrary toric variety the condition reformulates as
follows.
 The number of movable divisors (i. e. those that we can move from any point) is greater than codimension of
maximal cone whose edges are not subset of basis.


    \item[ii]
\label{remark:great_Picard_number} We assume that the Picard number
of the ambient variety is $1$ just for simplicity. In cases of
higher Picard number theorems~\ref{theorem:main}
and~\ref{theorem:toric} may be formulated and proved in a similar
manner. We consider a multivariable $q=(q_1,\ldots,q_m)$ instead of
the variable $q$ (here $m=\rk \pic X$) a multidegree of a curve
(with respect to the generators of Picard group) instead of a degree
and so on. The only difference in this case is the expression for
correction term $e^{-\alpha_X q}$.

    \item[iii]
\label{remark:dimension_3} If we lift the hypothesis $\pic X = \ZZ $
(or $\dim X\geq 3$ in remark~\ref{remark:Lefschetz_condition}), then
theorems~\ref{theorem:main} and~\ref{theorem:toric} give \emph{the
restricted} $I$-series (see
remark~\ref{remark:restricted_I-series}).

    \item[iiii]
The main theorems can be generalized to Calabi--Yau varieties. Their
proofs differ from proofs of theorems~\ref{theorem:main}
and~\ref{theorem:toric} only in formulas for correction term (see~\cite{Gi97}). 
\end{description}

\subsubsection{} {Examples.}
\begin{description}
    \item[i]
Consider a general hypersurface $X^k\subset \PP^n$. Let $k<n$. Then,
by theorem~\ref{theorem:main},
$$
    I^X=\sum_{d=0}^\infty q^d \frac{
    \os kH+1\cs_{\boldsymbol{kd+1}}
    }{
    (\os H+1\cs_{\boldsymbol{d+1}})^{n+1}
    }.
$$
This coincides with the formula for $I$-series for a hypersurface
in projective space from~\cite{Gi96}.

If $k=n$, then, by theorem~\ref{theorem:main},
$$
    I^X=e^{-k!q}\sum_{d=0}^\infty q^d \frac{
    \os kH+1\cs_{\boldsymbol{kd+1}}
    }{
    (\os H+1\cs_{\boldsymbol{d+1}})^{n+1}
    }.
$$
This coincides with the formula 4.2.2 from~\cite{Pr04}.

    \item[ii]
For projective space ($w_i=1$) theorem~\ref{theorem:main} is
Givental's Theorem for complete intersections in projective space.
\end{description}

\section{Appendix: Golyshev's conjecture}
\subsection{} Consider \emph{a family} of counting matrices $A^\lambda=A+\lambda
\mathbf{1}$, where $\mathbf{1}$ is the identity matrix
(see~\ref{definition:counting-matrix} for the definition of counting
matrix). Consider the one-dimensional torus $\mathbb G_m=
\mathrm{Spec}\, \CC [z, z^{-1}]$ and the differential operator
$D=z\frac{\partial}{\partial z}$. To construct the family of
matrices $M^\lambda$, put its elements $m_{kl}^\lambda$ as follows
$$
    m_{kl}^\lambda=
  \begin{cases}
    0, & \text{if $k>l+1$}, \\
    1, & \text{if $k=l+1$}, \\
    a_{kl}^\lambda\cdot (-\frac{\partial}{\partial z})^{l-k+1}, & \text{if $k<l+1$}.
  \end{cases}
$$
Consider the family of differential operators
$$
    \widetilde{L}^\lambda=\mathrm{det}_{\mathrm{right}}(D\mathbf{1}-M^\lambda),
$$
where $\mathrm{det}_{\mathrm{right}}$ means ``right determinant'',
i. e. the determinant, which is calculated with respect to \emph{the
rightmost} column; all minors are calculated in the same way. Divide
$\widetilde{L}^\lambda$ on the left by $D$. We get the family of
operators $L^\lambda$, so $\widetilde{L}^\lambda=DL^\lambda$.

\subsection{} \Def [\cite{Go02}, 1.8].  {\it The equation of the family $L^\lambda[\Phi(z)]=0$ is called
\emph{counting equation $D3$}. }

\subsection{} {\bf Golyshev's conjecture}
[V.\,Golyshev,~\cite{Go02}]. {\it The solutions of $D3$ equations
for Fano threefolds with Picard group $\ZZ$ are modular. More
precisely, let $X$ be such variety, $i_X$ be its index, and
$N=\frac{\deg X}{2i_X^2}$. Then in the family of counting equations
for $X$ there is one, $L^{\lambda_X}[\Phi(z)]=0$, whose solution is
an Eisenstein series of weight $2$ on $X_0(N)$. }

\medskip

More precise description of counting equations as Picard--Fuchs
equations see in~\cite{Go02}.

\subsection{} Based on this conjecture, Golyshev gives a list of predictions of counting matrices
of Fano threefolds. To check this conjecture one should find all
such matrices.

There are $17$ families of smooth Fano threefolds with Picard
group
$\ZZ$ (V.\,Iskovskikh, \cite{Is77}, \cite{Is78}). 
For 14 of them counting matrices were found by
A.\,Beauville~(\cite{Bea95}), A.\,Kuznetsov (unpublished),
V.\,Golyshev (unpublished) and the author~(\cite{Pr04}). The
matrices for 2 other varieties are given by
theorem~\ref{theorem:matrices}.

All these matrices coincide with ones that were predicted by
Golyshev. Thus, theorem~\ref{theorem:matrices} finishes the proof
of Golyshev's conjecture.

\medskip

The author is grateful to I.\,Cheltsov, S.\,Galkin, A.\, Givental,
V.\,Golyshev, M.\,Kazarian, S.\,Shadrin, K.\,Shramov,
M.\,Tsfasman, and F.\,Zak for comments.

\end{document}